\documentclass[12pt]{amsart}

\usepackage[all]{xypic}
\usepackage{tikz}
\usepackage{ragged2e}
\usetikzlibrary{arrows} 
\usetikzlibrary{decorations.markings}
\usepackage{graphicx}
\usepackage{bm}
\usepackage{float}
\usepackage{epsf}
\usepackage{verbatim} 
\usepackage{amsmath}
\usepackage{amsfonts}
\usepackage{amssymb}
\usepackage{mathrsfs}
\usepackage{amsthm}
\usepackage{newlfont}
\usepackage{enumitem}
\usepackage[new]{old-arrows}
\usepackage{booktabs}
\usepackage{enumitem}
\usepackage{makecell}
\usepackage[multiple]{footmisc}
\usepackage[answerdelayed,lastexercise]{exercise}
\usepackage[hidelinks]{hyperref}
\usepackage{mleftright}

\usepackage{fnpct}
\usepackage{colonequals} 
\usepackage{stmaryrd}

\usepackage{apptools}
\usepackage{chngcntr}

\usepackage[margin=1in]{geometry}

\newtheorem{thm}{Theorem}[section]
\newtheorem{prop}[thm]{Proposition}
\newtheorem{lem}[thm]{Lemma}
\newtheorem{cor}[thm]{Corollary}

\newtheorem{thmx}{Theorem}

\theoremstyle{definition}

\newtheorem{example}[thm]{Example}

\theoremstyle{remark}
\newtheorem{rem}[thm]{Remark}

\numberwithin{equation}{section}
\numberwithin{thm}{section}

\makeatletter 
\newcommand\mynobreakpar{\vspace{0.02in}\par\nobreak\@afterheading}  
\makeatother

\makeatletter
\@addtoreset{footnote}{section}
\makeatother

\usepackage{xcolor}
\hypersetup{
	colorlinks,
	linkcolor={red!50!black},
	citecolor={blue!50!black},
	urlcolor={blue!80!black}
}

\AtBeginDocument{%
	\def\MR#1{}
}

\DeclareMathOperator{\Hom}{Hom}
\DeclareMathOperator{\End}{End}

\DeclareMathOperator{\rank}{rank}

\DeclareMathOperator{\inv}{inv}

\DeclareMathOperator{\chr}{char}

\DeclareMathOperator{\Mat}{Mat}
\DeclareMathOperator{\ord}{ord}

\DeclareMathOperator{\Gal}{Gal}

\DeclareMathOperator{\Nr}{Nr}

\DeclareMathOperator{\Cent}{Cent}

\newcommand{\sep}{\mathrm{sep}}

\newcommand{\fP}{\mathfrak{P}}
\newcommand{\fQ}{\mathfrak{Q}}

\newcommand{\fl}{\mathfrak{l}}

\newcommand{\fp}{\mathfrak{p}}
\newcommand{\fq}{\mathfrak{q}}

\newcommand{\cM}{\mathcal{M}}

\newcommand{\cO}{\mathcal{O}}

\newcommand{\F}{\mathbb{F}}

\newcommand{\Q}{\mathbb{Q}}

\newcommand{\Z}{\mathbb{Z}}

\newcommand{\To}{\longrightarrow}

\newcommand{\Fi}{F_\infty}

\newcommand{\oF}{\overline{\F}}
\newcommand{\oK}{\overline{K}}

\newcommand{\La}{\Lambda}
\newcommand{\la}{\lambda}

\newcommand{\twist}[1]{#1\!\left\{\tau\right\}}

\newcommand{\twistdvr}[1]{#1\!\left\{\mkern-4mu\left\{\tau\right\}\mkern-4mu\right\}}

\newcommand{\ls}[2]{#1\!\left(\mkern-4mu\left(#2\right)\mkern-4mu\right)}
\newcommand{\dvr}[2]{#1\!\left\llbracket #2\right\rrbracket}

\newcommand{\Mod}[1]{\ (\mathrm{mod}\ #1)}
\newcommand{\abs}[1]{\left|#1\right|}

\newcommand{\prodprime}{\sideset{}{'}\prod}

\theoremstyle{plain}
\newtheorem{thmxp}{Theorem}

\title{Endomorphism algebras of Drinfeld modules in special characteristic}

\author{Mihran Papikian}
\address{Department of Mathematics, Pennsylvania State University, University Park, Pennsylvania, United States of America}
\email{papikian@psu.edu}

\thanks{The author was supported in part by the Simons Foundation, award number MPS-TSM-00008093.}


\begin{document}
\begin{abstract}
	We classify the division algebras that occur as endomorphism algebras of Drinfeld modules of a given rank in special 
	characteristic which cannot be defined over a finite field. The classification is in terms of the center of the algebra 
	and its local invariants. Together with the known results in generic characteristic and over finite fields, this gives a 
	complete description of the endomorphism algebras of Drinfeld modules of a given rank. 
\end{abstract}
\maketitle

\section{Introduction}

\subsection{Background} The endomorphism algebra of a simple abelian variety is a division algebra over $\Q$ with a positive
involution. Such algebras were classified by Albert, and the question of which of them occur for abelian varieties of a given
dimension is a classical one, with many applications to the arithmetic and the moduli of abelian varieties. In characteristic
zero it was answered by Albert and Shimura, and over finite fields by Honda and Tate. Over arbitrary fields of positive
characteristic every such algebra occurs, but the dimensions in which a given algebra occurs are known only in special cases;
see \cite{OortEnd}, \cite[(3.5)]{ChaiOort}. This paper concerns the analogous question for Drinfeld modules, where a complete
answer is possible.

Let $F$ be a global function field with constant field $\F_q$, let $\infty$ be a place of $F$, and let
$A\subset F$ be the ring of elements regular away from $\infty$. Let $K$ be an $A$-field, i.e., a field equipped with an 
$\F_q$-algebra homomorphism $\gamma\colon A\to K$, and let $\phi$ be a Drinfeld $A$-module of rank $r$ over $K$. The endomorphism 
algebra $\End^\circ_K(\phi)=\End_K(\phi)\otimes_A F$ is a division algebra of finite dimension over $F$ 
\cite[Cor. to Prop. 2.4]{Drinfeld74}. The possible structures of this algebra depend on the $A$-characteristic $\chr_A(K)=\ker(\gamma)$.

If $\chr_A(K)=0$, then $\End^\circ_K(\phi)$ is a field extension of $F$ in which $\infty$ does not split \cite[Cor. to Prop. 2.4]{Drinfeld74}, 
and its degree divides $r$ (see Proposition \ref{propInvInf} below). Every such field occurs, as one sees from lattices over orders 
of the field and the analytic uniformization \cite[Prop. 3.1]{Drinfeld74}, \cite[Prop. 2.4]{Gekeler83}; for $A=\F_q[T]$ 
see \cite[$\S$5.3]{PapikianGTM}.

If $K$ is finite, then $\End^\circ_K(\phi)$ is a central division algebra over the field $F(\pi)$ generated by the Frobenius endomorphism $\pi$.
Its local invariants are determined by $r$ and the two places of $F(\pi)$ where $\pi$ has nonzero valuation.
Together with the analogue of the Honda--Tate theorem, this gives a complete description of the algebras that occur
\cite{Drinfeld77, GekelerFDM, Yu95}, \cite[Ch. 2]{Laumon}.

The remaining case, $\chr_A(K)=\fp\neq 0$ with $\phi$ not definable over a finite field, lies between these two.
Here $\End^\circ_K(\phi)$ can be noncommutative \cite[Ex. 5.2]{PinkII}, but by a theorem of Pink \cite[Thm. 6.1]{PinkII},
a maximal subfield cannot have degree $r$ over $F$. Pink's work \cite{PinkI, PinkII} and that of Devic and Pink \cite{DevicPink}
describe the image of Galois on the Tate modules of $\phi$ in terms of the endomorphism rings of $\phi$ and of its
restrictions to subrings of $A$. See also \cite{PinkKummer, KuhnPink}, where the structure of these rings is used and computed. 
The question of which algebras occur as $\End^\circ_K(\phi)$ in this case arose during the writing of \cite{PapikianGTM}, 
where it was noted that a complete classification seemed to be lacking \cite[Rem. 4.1.5 (5)]{PapikianGTM}. To the author's 
knowledge, it has remained open. The purpose of this paper is to fill this gap.

\subsection{Main results} We say that $\phi$ is \textit{isotrivial} if it becomes isomorphic, over some field extension of $K$,
to a Drinfeld module defined over a finite field. We put $\End^\circ(\phi)\colonequals \End^\circ_{\oK}(\phi)$. By a
\textit{division $F$-algebra} we mean a division ring $D$ of finite dimension over a subfield $F$ of its center. For such $D$,
let $Z$ be the center of $D$ and put
$$
e(D)\colonequals [Z:F], \qquad d(D)\colonequals \sqrt{[D:Z]}, 
$$
so that 
$$
[D:F]=e(D)\, d(D)^2.
$$

\begin{thmx}\label{thmA}
	Let $\fp\lhd A$ be a nonzero prime, let $r\geq 1$, and let $D$ be a division $F$-algebra with center $Z$. Put $e=e(D)$ and $d=d(D)$.
	There are an $A$-field $K$ with $\chr_A(K)=\fp$ and a non-isotrivial Drinfeld module $\phi$ of rank $r$ over $K$
	such that $\End^\circ(\phi)\cong D$ as $F$-algebras if and only if the following hold:
	\begin{enumerate}
		\item $ed$ divides $r$, and $[D:F]<r$;
		\item $Z$ has a unique place $\infty_Z$ over $\infty$;
		\item there is a place $\fp_Z$ of $Z$ over $\fp$ such that $D$ splits at every place of $Z$ other than $\infty_Z$ and $\fp_Z$;
		\item $s\cdot\inv_{\infty_Z}(D)=-1/d$ in $\Q/\Z$, where $s=r/(ed)$.
	\end{enumerate}
	In this case $K$ can be taken finitely generated over $\F_q$.
\end{thmx}

If (1)--(4) hold, then $\gcd(d,s)=1$, the invariant $\inv_{\infty_Z}(D)\in\frac1d\Z/\Z$ is determined by $s$, and
$\inv_{\fp_Z}(D)=-\inv_{\infty_Z}(D)$. So $D$ is determined up to isomorphism
by $Z$, $\fp_Z$, $d$ and $r$. Conversely, given $Z$ with a unique place over $\infty$, a place $\fp_Z$ of $Z$ over $\fp$, and 
$d\geq1$ with $ed\mid r$, $ed^2<r$ and $\gcd(d,r/(ed))=1$, the central division algebra over $Z$ with these invariants occurs. 

For $d=1$, Theorem \ref{thmA} says that a field $Z$ occurs if and only if $\infty$ does not split in $Z$ and $[Z:F]$ is a proper 
divisor of $r$. There is no condition at $\fp$. In generic characteristic the same fields occur, together with those of degree $r$. 
Condition (1) sharpens Pink's theorem \cite[Thm. 6.1]{PinkII}, which says that $ed<r$;
see Corollary \ref{corPink}. For example, Theorem \ref{thmA} (1) implies that the endomorphism algebra of a non-isotrivial
Drinfeld module of rank $\leq5$ is commutative, whereas quaternion algebras occur in rank $6$.

In terms of the \textit{endomorphism ratio} $\mu(\phi)\colonequals r/[\End^\circ(\phi):F]$, Theorem \ref{thmA} implies that
$\mu(\phi)>1$ for non-isotrivial $\phi$, that every rational number $>1$ occurs, and that $\mu(\phi)\leq 1$ if and only if
$\phi$ is isotrivial; see Corollary \ref{corRatio}. This answers the analogue for Drinfeld modules of a question of Oort on
abelian varieties \cite[Question 22.2]{OortQ}. For abelian varieties in characteristic $p$ the answer is different: by
\cite{Tafazolian}, every positive rational number occurs.

Theorem \ref{thmA} is deduced from the following criterion. 

\begin{thmxp}\label{thmAp}
	Let $\fp\lhd A$ be a nonzero prime, let $D$ be a division $F$-algebra, and let $r\geq 1$. The following are equivalent:
	\begin{enumerate}
		\item There are an $A$-field $K$ with $\chr_A(K)=\fp$ and a non-isotrivial Drinfeld module $\phi$ of rank $r$ over $K$
		such that $\End^\circ(\phi)\cong D$ as $F$-algebras.
		\item There are a Drinfeld module $\psi$ of rank $r'<r$ over $\oF_\fp$ and an embedding of $F$-algebras
		$D\hookrightarrow\End^\circ(\psi)$ such that $[D:F]$ divides $r-r'$.
	\end{enumerate}
	In (1) the field $K$ can be taken finitely generated over $\F_q$, and in (2) one can take $r'=r-[D:F]$. Moreover, if (1) holds, 
	then every $\psi$ of rank $r'$ over $\oF_\fp$ with $D\hookrightarrow\End^\circ(\psi)$ satisfies $r\equiv r'\Mod{[D:F]}$.
\end{thmxp}

To pass from Theorem \ref{thmAp} to Theorem \ref{thmA}, we describe the division subalgebras of the endomorphism algebra of a 
Drinfeld module over a finite field in terms of local invariants (Proposition \ref{propFF}). Combined with the analogue of the 
Honda--Tate theorem, this shows that $D$ embeds into $\End^\circ(\psi)$ for some $\psi$ of rank $r'$ over $\oF_\fp$ if and only if 
$ed$ divides $r'$ and conditions (2)--(4) of Theorem \ref{thmA} hold with $r'$ in place of $r$ (Corollary \ref{corFFinv}). The proof of Theorem \ref{thmAp} gives the same statement for endomorphism rings:
an $A$-algebra $\cO$ is the endomorphism ring over $\oK$ of a non-isotrivial Drinfeld module of rank $r$ and characteristic $\fp$
if and only if $\cO$ embeds as an $A$-subalgebra into $\End_{\oF_\fp}(\psi)$ for some $\psi$ over $\oF_\fp$ of rank $r'<r$ with
$r\equiv r'$ modulo $\dim_F(\cO\otimes_AF)$; see Theorem \ref{thmRings}. Endomorphism rings are not determined by local 
invariants alone, and we do not attempt an analogue of Theorem \ref{thmA} for them. 

The last assertion of Theorem \ref{thmAp} and the necessity of conditions (2)--(4) in Theorem \ref{thmA} are special cases of the following result, which holds in every characteristic.

\begin{thmx}\label{thmB}
	Let $\phi$ be a Drinfeld module of rank $r$ over an arbitrary $A$-field $K$, and let $D\subseteq\End^\circ(\phi)$ be a division
	$F$-subalgebra with center $Z$. Put $e=e(D)$, $d=d(D)$ and $s=r/(ed)$.
	\begin{enumerate}
		\item $s$ is an integer, $Z$ has a unique place $\infty_Z$ over $\infty$, and $D\otimes_Z Z_{\infty_Z}$ is a division algebra
		whose invariant satisfies
		$$
		s\cdot\inv_{\infty_Z}(D)=-\frac{1}{d}\quad\text{in }\Q/\Z.
		$$
		In particular $\gcd(d,s)=1$.
		\item There is at most one place of $Z$ other than $\infty_Z$ at which $D$ does not split, and such a place lies over $\chr_A(K)$.
		\item If $D$ also embeds into $\End^\circ(\phi')$ for a Drinfeld module $\phi'$ of rank $r'$ over some $A$-field,
		then $r\equiv r'\Mod{[D:F]}$.
	\end{enumerate}
\end{thmx}

When $\phi$ is defined over a finite field and $D=\End^\circ(\phi)$, we have $s=1$, and (1) recovers the known value
$\inv_{\infty_Z}(D)=-1/d$. For non-isotrivial $\phi$ the invariant at $\infty$ can be different from $-1/d$.
By Theorem \ref{thmAp}, $\End^\circ(\phi)$ embeds into the endomorphism algebra of a Drinfeld module over a finite field. 
One might expect that $\End^\circ(\phi)$ is itself the endomorphism algebra of a Drinfeld module over a finite field of 
the same $A$-characteristic. The following result shows that this is not the case, even for commutative algebras.

\begin{thmx}\label{thmC}
	Let $A=\F_q[T]$, so that $F=\F_q(T)$ and $1/T$ is a uniformizer at $\infty$. There are non-isotrivial Drinfeld modules 
	$\phi$ of $A$-characteristic $(T)$ whose endomorphism algebras are not endomorphism algebras of Drinfeld modules over 
	finite fields of $A$-characteristic $(T)$. More precisely:
	\begin{enumerate}
		\item For odd $q$ and a non-square $c\in\F_q^\times$, there is an explicit $\phi$ of rank $4$ over $\F_{q^2}(x)$ with
		$\End^\circ(\phi)\cong F(\sqrt{T+c})$. This field is not the endomorphism algebra of a Drinfeld module over a finite 
		field of $A$-characteristic $(T)$, but it is the endomorphism algebra of a Drinfeld module over a finite field of 
		$A$-characteristic $\fq$ for every prime $\fq$ of $A$ which splits in it.
		\item There is an explicit $\phi$ of rank $15$ over $\F_{q^3}(x)$ such that $\End^\circ(\phi)$ is the central division algebra
		over $F$ of dimension $9$ with invariants $-1/3$ at $(T)$ and $1/3$ at $\infty$. This algebra is not the endomorphism
		algebra of a Drinfeld module over a finite field of any $A$-characteristic.
	\end{enumerate}
\end{thmx}

In (1), $\infty$ has a unique place in $Z=F(\sqrt{T+c})$, so does $(T)$, and $Z$ is not a constant field extension of $F$. 
The first property is needed for $Z$ to occur at all, by Theorem \ref{thmA} (2). The other two prevent $Z$ from being the 
endomorphism algebra of a Drinfeld module over a finite field of $A$-characteristic $(T)$; see Example \ref{ex1} (3). 
In (2), the algebra is the opposite of the endomorphism algebra of the supersingular Drinfeld module $\psi$ with 
$\psi_T=\tau^3$. Its invariant at $\infty$ is $1/3$, whereas over a finite field a Drinfeld module with endomorphism 
algebra of center $F$ and dimension $9$ has invariant $-1/3$ at $\infty$; see Example \ref{ex2} (1).

\subsection{Methods} The proofs use only standard tools from the theory of Drinfeld modules and central simple algebras.
Theorem \ref{thmB} (1) comes from the embedding of $\End(\phi)\otimes_A\Fi$ into the central division algebra of
invariant $-1/r$ over $\Fi$ formed by twisted Laurent series in $\tau^{-1}$. 
For the direction (1) $\Rightarrow$ (2) of Theorem \ref{thmAp}, we degenerate $\phi$ at a discrete valuation where it does not 
have potentially good reduction, and use the injectivity of reduction on endomorphisms \cite[Lem. 3.3]{Gekeler83}. Then we induct 
on the transcendence degree of $K$. The congruence follows from Theorem \ref{thmB}. For (2) $\Rightarrow$ (1), we construct $\phi$ 
by Tate uniformization over $\ls{k}{x}$ from $\psi$ and a lattice which is a projective module over an order of $D$. Here we use 
Drinfeld's Tate uniformization \cite[$\S$7]{Drinfeld74} and Gekeler's description of the endomorphisms of a Tate-uniformized module 
\cite[Prop. 3.8]{Gekeler83}, both of which are valid in any characteristic. The same degeneration technique appears, for Drinfeld 
modules of rank $2$ over $\F_q[T]$, in the proof of \cite[Prop. 2.12]{KuhnPink}. For abelian varieties, the analogue of
the construction in (2) $\Rightarrow$ (1), with Tate uniformization replaced by the Mumford--Faltings construction of
degenerating abelian varieties, is due to Oort and van der Put \cite[Thm. (1.1)]{OortvdPut}. The relation $r=r'+n[D:F]$
corresponds to $\dim(Y)=\dim(Z)+m[D:\Q]$ there. The direction (1) $\Rightarrow$ (2) has no counterpart for abelian varieties,
since a non-isotrivial abelian variety can have potentially good reduction everywhere, whereas a non-isotrivial Drinfeld
module over a finitely generated field cannot (Lemma \ref{lemIsotriv} (4)).

Theorem \ref{thmA} then follows from
Theorem \ref{thmAp}, Theorem \ref{thmB}, the Honda--Tate theorem for Drinfeld modules \cite[Thm. 3]{Yu95}, \cite[(2.2.2)]{Laumon}, 
and the local criterion for a field to embed into a central division algebra over a global field \cite[(A.3.4)]{Laumon}. 
Theorem \ref{thmB} (2) is a byproduct. 
None of the results of \cite{PinkI, PinkII, DevicPink} on the Galois representations is used.

\subsection{Notation and conventions}\label{ssNotation} Throughout, $F$, $\infty$ and $A$ are as above, $\Fi$ is the completion of $F$ at $\infty$, 
and $\abs{\cdot}$ is the absolute value on $\Fi$ normalized so that $\abs{a}=\#(A/aA)$ for $0\neq a\in A$; cf. \cite[(1.0)]{Laumon}. 
Put $\deg(a)=\log_q\abs{a}$. For an $A$-field $K$, let $\twist{K}$ be the twisted polynomial ring with $\tau c=c^q\tau$, and let 
$\partial\colon\twist{K}\to K$ be the constant term. A Drinfeld module of rank $r$ over $K$ is an $\F_q$-algebra homomorphism 
$\phi\colon A\to\twist{K}$, $a\mapsto\phi_a$, such that $\partial\phi_a=\gamma(a)$ and $\deg_\tau\phi_a=r\deg(a)$ for all 
$a\in A$ \cite[$\S$2]{Drinfeld74}, \cite[(1.2.1)]{Laumon}. 

For a field extension $L/K$, we write $\End_L(\phi)$ for the endomorphisms defined over $L$. If $L\supseteq K^\sep$, then
$\End_L(\phi)=\End_{K^\sep}(\phi)$. This is \cite[Thm. 3.3.9]{PapikianGTM}, whose proof uses only an element $a\in A$ prime to 
$\chr_A(K)$ with $\#\phi[a]$ large, and so applies to any $A$. The ring $\End_{K^\sep}(\phi)$ is a projective $A$-module of rank 
at most $r^2$, and $\End_{K^\sep}(\phi)\otimes_A\Fi$ is a division algebra \cite[Prop. 2.4 and Cor.]{Drinfeld74}, 
\cite[(2.1.4)]{Laumon}. Hence $\End_{K^\sep}(\phi)=\End_{K'}(\phi)$ for some finite separable extension $K'/K$. In particular, 
$\End^\circ(\phi)$ does not change if $K$ is replaced by an overfield, or by a subfield over which $\phi$ is defined.

Local invariants of central simple algebras over local fields take values in $\Q/\Z$ and are normalized as in 
\cite[(A.2.2)]{Laumon} and \cite[(2.1)]{GekelerFDM}. With this normalization, if $\phi$ is defined over a finite 
field $k$ with Frobenius $\pi$, then $\End^\circ_k(\phi)$ is a central division algebra over $F(\pi)$ of degree $r/[F(\pi):F]$, the 
field $F(\pi)$ has a unique place over $\infty$ and a unique place dividing $\pi$, and the invariants of $\End^\circ_k(\phi)$ are 
$-[F(\pi):F]/r$ at the first of these places, $[F(\pi):F]/r$ at the second, and $0$ elsewhere \cite[Thm. 2.9]{GekelerFDM}, 
\cite[Thm. 1]{Yu95}, \cite[(2.2.2)]{Laumon}.

The facts about central simple algebras that we use are collected in $\S$\ref{sCSA}. 


\section{Central simple algebras}\label{sCSA}

We recall the facts about central simple algebras that are used below. The standard reference is \cite[$\S\S$7, 31, 32]{Reiner}. 
The statements in the form we need are in \cite[App. A]{Laumon}. Throughout this section, $k$ is a field and $\Delta$ is a central 
simple $k$-algebra.

\begin{itemize}
	\item[(C1)] \textit{Double centralizer.} Let $B\subseteq\Delta$ be a simple $k$-subalgebra. Then $\Cent_\Delta(B)$ is simple, and 
	$$
	\Cent_\Delta(\Cent_\Delta(B))=B,\qquad [\Delta:k]=[B:k]\cdot[\Cent_\Delta(B):k]
	$$
	\cite[(A.1.2)]{Laumon}. If the center of $B$ is $k$, then $B\otimes_k\Cent_\Delta(B)$ is simple, so the multiplication map to 
	$\Delta$ is injective, and by comparing dimensions 
	$$
	\Delta\cong B\otimes_k\Cent_\Delta(B).
	$$
	\item[(C2)] \textit{Subfields.} Let $\Delta$ be a division algebra of degree $m$, and let $L\subseteq\Delta$ be a subfield 
	containing $k$, with $\ell=[L:k]$. Then $\ell$ divides $m$, and $\Cent_\Delta(L)$ is a central division algebra over $L$ of degree 
	$m/\ell$ such that 
	$$
	\Delta\otimes_kL\cong\Mat_\ell(\Cent_\Delta(L))
	$$
	\cite[(A.1.4)]{Laumon}. 
	\item[(C3)] \textit{Local invariants.} Let $k$ be a local field. The invariant is an isomorphism $\inv\colon\mathrm{Br}(k)\to\Q/\Z$, 
	and the invariant of a central division algebra of degree $n$ has order $n$ \cite[(A.2.2)]{Laumon}. In the situation of (C2), 
	$$
	\inv(\Cent_\Delta(L))=[L:k]\cdot\inv(\Delta)
	$$
	\cite[(A.2.4)]{Laumon}. 
	\item[(C4)] \textit{Global invariants.} Let $k$ be a global function field, and for a place $v$ of $k$ put 
	$\inv_v(\Delta)=\inv(\Delta\otimes_kk_v)$. Then $\inv_v(\Delta)=0$ for all but finitely many $v$, $\sum_v\inv_v(\Delta)=0$, the Brauer class of 
	$\Delta$ is determined by the $\inv_v(\Delta)$, every family of invariants with these two properties occurs, and the index of 
	$\Delta$ is the least common multiple of the orders of the $\inv_v(\Delta)$ \cite[(A.3.3)]{Laumon}. 
	\item[(C5)] \textit{Restriction.} Let $k$ be a global function field, $L/k$ a finite extension, and $w$ a place of $L$ over a place $v$ of $k$. Then 
	$$
	\inv_w(\Delta\otimes_kL)=[L_w:k_v]\cdot\inv_v(\Delta)
	$$
	\cite[(A.3.2)]{Laumon}. 
	\item[(C6)] \textit{Embedding of fields.} Let $k$ be a global function field, $\Delta$ a division algebra of degree $m$, and $L/k$ a finite 
	extension of degree $\ell$. Then $L$ is $k$-isomorphic to a subfield of $\Delta$ if and only if $\ell$ divides $m$ and 
	$$
	\frac{m}{\ell}\,[L_w:k_v]\cdot\inv_v(\Delta)=0\qquad\text{for all places $w\mid v$}
	$$
	\cite[(A.3.4)]{Laumon}. 
\end{itemize}


\section{The invariant at $\infty$}\label{sInf}

In this section $K$ is an arbitrary $A$-field.

\begin{prop}\label{propInvInf}
	Let $\phi$ be a Drinfeld module of rank $r$ over $K$, and let $D\subseteq\End^\circ(\phi)$ be a division $F$-subalgebra
	with center $Z$. Put $e=e(D)$ and $d=d(D)$. Then:
	\begin{enumerate}
		\item $ed$ divides $r$. Put $s\colonequals r/(ed)$.
		\item $Z_\infty\colonequals Z\otimes_F\Fi$ is a field, and $D_\infty\colonequals D\otimes_Z Z_\infty$ is a central division algebra over $Z_\infty$.
		\item $s\cdot\inv(D_\infty)=-1/d$ in $\Q/\Z$. In particular, $\gcd(d,s)=1$.
	\end{enumerate}
\end{prop}
\begin{proof}
	We may assume that $K$ is algebraically closed. The homomorphism $\phi$ extends by continuity to an injection of $\Fi$ into the 
	ring $\ls{K}{\tau^{-1}}$ of twisted Laurent series in $\tau^{-1}$ \cite[$\S$3.1]{YuST}. Let $\Delta_\phi$ be the centralizer of 
	$\phi(\Fi)$ in $\ls{K}{\tau^{-1}}$. It is a central division algebra over $\Fi$ of degree $r$ with invariant $-1/r$ 
	\cite[$\S$1.4 and Lem. 2.4]{YuST}, \cite[Lem. 2.2]{Zywina}; for $A=\F_q[T]$ see also \cite[Prop. 3.4.12, Rem. 3.4.16 (3)]{PapikianGTM}. 
	Every $u\in\End(\phi)$ commutes with $\phi(A)$, hence with $\phi(\Fi)$, so $\alpha\otimes u\mapsto\phi(\alpha)u$ defines a 
	homomorphism of $\Fi$-algebras $\Fi\otimes_A\End(\phi)\to\Delta_\phi$. It is injective because $\Fi\otimes_A\End(\phi)$ is a 
	division algebra. Since $\Fi$ is flat over $F$, we obtain
	\begin{equation}\label{eqEmbInf}
	D\otimes_F\Fi\longhookrightarrow\End^\circ(\phi)\otimes_F\Fi\cong\Fi\otimes_A\End(\phi)\longhookrightarrow\Delta_\phi. 
	\end{equation}

	(1) Let $M$ be a maximal subfield of $D$, so $[M:F]=ed$. By \eqref{eqEmbInf}, $M\otimes_F\Fi$ is a commutative subalgebra of 
	the division algebra $\Delta_\phi$, hence a field of degree $ed$ over $\Fi$. By (C2), $ed$ divides $r$. 

	(2) By \eqref{eqEmbInf}, $Z_\infty$ is a subring of a division algebra, hence a field, and $[Z_\infty:\Fi]=e$. 
	Let $\Delta'\colonequals\Cent_{\Delta_\phi}(Z_\infty)$. This is a central division algebra over $Z_\infty$
	of degree $r/e=ds$ with invariant $e\cdot(-1/r)=-1/(ds)$, by (C2) and (C3). The algebra $D_\infty=D\otimes_F\Fi$ is central simple over
	$Z_\infty$ of degree $d$ and is contained in $\Delta'$, so it is a division algebra.

	(3) Let $C\colonequals\Cent_{\Delta'}(D_\infty)$. By (C1), $C$ is a central division algebra over $Z_\infty$ of degree $s$, and
	$\Delta'\cong D_\infty\otimes_{Z_\infty}C$. By (C3), 
	$$
	\inv(D_\infty)+\inv(C)=-\frac{1}{ds}\qquad\text{and}\qquad s\cdot\inv(C)=0.
	$$
	Multiplying the first equality by $s$ gives $s\cdot\inv(D_\infty)=-1/d$. Since $\inv(D_\infty)$ lies in $\frac1d\Z/\Z$,
	this shows that the generator $-1/d$ of this cyclic group of order $d$ is a multiple of $s$, so $\gcd(d,s)=1$.
\end{proof}

\begin{cor}\label{corRankCong}
	Let $D$ be a division $F$-algebra. Suppose that $D$ embeds as an $F$-subalgebra into $\End^\circ(\phi_1)$ and into $\End^\circ(\phi_2)$,
	where $\phi_i$ is a Drinfeld module of rank $r_i$ over an $A$-field $K_i$ for $i=1,2$. Then $r_1\equiv r_2\Mod{[D:F]}$.
\end{cor}
\begin{proof}
	The invariant $\iota\colonequals\inv(D_\infty)$ of Proposition \ref{propInvInf} depends only on $D$, and it has order $d$ 
	by (C3), because $D_\infty$ is a division algebra of degree $d$. Let $s_i=r_i/(ed)$. Then $s_1\iota=-1/d=s_2\iota$, so 
	$(s_1-s_2)\iota=0$ and $d$ divides $s_1-s_2$. Hence $r_1-r_2=ed(s_1-s_2)$ is divisible by $ed^2=[D:F]$.
\end{proof}

\begin{rem}\label{remInvInf}
	\begin{enumerate}
		\item Suppose that $\phi$ is defined over a finite field $k'$ and $D=\End^\circ_{k'}(\phi)$. Then $r=ed$ \cite[Thm. 2.9]{GekelerFDM}, 
		so $s=1$ and Proposition \ref{propInvInf} (3) gives $\inv(D_\infty)=-1/d$, in agreement with \cite[Thm. 2.9]{GekelerFDM}. 
		Thus a division algebra whose invariant at $\infty$ is not $-1/d$ cannot be the endomorphism algebra of a Drinfeld module 
		over a finite field. This is the source of the second example in Theorem \ref{thmC}.
		\item In generic characteristic $D$ is a field, $d=1$, and Proposition \ref{propInvInf} only says that $\infty$ does not split in $Z$
		and that $[Z:F]$ divides $r$. The algebra $\Delta'$ in the proof is then the central division algebra of invariant $-e/r$ over 
		$Z_\infty$ that appears in the Sato--Tate law for Drinfeld modules: it is the algebra $B_\phi$ of \cite[$\S$1.2]{Zywina}. 
		Yu \cite{YuST} works instead with $\Delta_\phi$ itself. 
	\end{enumerate}
\end{rem}

\begin{proof}[Proof of Theorem \ref{thmB} (1) and (3)]
	Part (1) is Proposition \ref{propInvInf}, and (3) is Corollary \ref{corRankCong}. 
\end{proof}

Theorem \ref{thmB} (2) is vacuous in generic characteristic, where $D$ is commutative. In special characteristic it is proved 
in $\S$\ref{sFF}; see Corollary \ref{corSplit}. 


\section{Isotriviality and stable reduction}\label{sIsotriv}

From now on, $K$ is an $A$-field with $\chr_A(K)=\fp\neq0$. Then $\gamma(A)\cong A/\fp$ is a finite field. 
Fix a nonconstant $a\in A$ and put $N=r\deg(a)$. Write
$$
\phi_a=\gamma(a)+g_1\tau+\cdots+g_N\tau^N,\qquad g_N\neq0,
$$
and put  
$$
J_i(\phi)\colonequals\frac{g_i^{q^N-1}}{g_N^{q^i-1}}\in K\qquad (1\leq i\leq N-1).
$$
For $A=\F_q[T]$ and $a=T$, these are $J$-invariants in the sense of Potemine \cite[Def. 2.1]{Potemine}: $J_i(\phi)$ is the 
$(q^{\gcd(i,r)}-1)$-th power of his $j_i$-invariant \cite[(2.5)]{Potemine}; see also \cite[\S3.8]{PapikianGTM}.

Let $v$ be a discrete valuation on $K$, with valuation ring $R$, maximal ideal $\cM$ and residue field $k$. Note that $\gamma(A)\subset R$.
We say that $\phi$ has \textit{stable reduction} (resp. \textit{good reduction}) at $v$ if there is $c\in K^\times$ such that all 
$c^{-1}\phi_bc$, $b\in A$, have coefficients in $R$ and their reductions modulo $\cM$ form a Drinfeld module (resp. a Drinfeld 
module of rank $r$) over $k$ \cite[$\S$3]{Gekeler83}. We say that $\phi$ has \textit{potentially good reduction} at $v$ if it has 
good reduction over some finite extension of $K$, with respect to some extension of $v$. Put
\begin{align*}
e(\phi) &=\min_{1\leq i\leq N}\frac{v(g_i)}{q^i-1},\\  
N_v(\phi) &=\max\left\{1\leq i\leq N\ \Big|\ \frac{v(g_i)}{q^i-1}=e(\phi)\right\}, \\ 
r_v(\phi) &=\frac{N_v(\phi)}{\deg(a)}.
\end{align*}

\begin{lem}\label{lemStable}
	Let $L/K$ be a finite extension and $w$ an extension of $v$ to $L$. Then $\phi$ has stable reduction over $L$ at $w$ if and only if 
	$e(\phi)\in w(L^\times)$. Every stable reduction of $\phi$, over any finite extension of $K$ and with respect to any extension of $v$, 
	has rank $r_v(\phi)$. In particular $r_v(\phi)$ is an integer, and $\phi$ has potentially good reduction at $v$ if and only if $r_v(\phi)=r$.
\end{lem}
\begin{proof}
	Let $R_w$ and $k_w$ be the valuation ring and the residue field of $w$. Let $c\in L^\times$ and put $\psi=c^{-1}\phi c$, so that $\psi_a=\gamma(a)+\sum_ic^{q^i-1}g_i\tau^i$ and 
	$$
	w(c^{q^i-1}g_i)=(q^i-1)\Bigl(w(c)+\frac{v(g_i)}{q^i-1}\Bigr). 
	$$
	Hence $\psi_a$ has coefficients in $R_w$ and nonconstant reduction if and only if $w(c)=-e(\phi)$, and in this case the reduction of 
	$\psi_a$ has $\tau$-degree $N_v(\phi)$. If $\psi$ has coefficients in $R_w$ and its reduction $\bar\psi$ is a Drinfeld module, then $\bar\psi_a$ is nonconstant, so 
	$w(c)=-e(\phi)$ and $\bar\psi$ has rank $N_v(\phi)/\deg(a)=r_v(\phi)$. 

	Conversely, let $w(c)=-e(\phi)$. For $b\in A$, we have $\psi_b\psi_a=\psi_a\psi_b$, where $\psi_a\in\twist{R_w}$ has nonconstant 
	reduction, so $\psi_b\in\twist{R_w}$ by \cite[Lem. 3.3 (i)]{Gekeler83}. The reductions $\bar\psi_b$ form an $\F_q$-algebra 
	homomorphism $\bar\psi\colon A\to\twist{k_w}$ with $\partial\bar\psi_b=\overline{\gamma(b)}$ and $\bar\psi_a$ nonconstant. 
	Such a homomorphism is a Drinfeld module \cite[Prop. 2.1 and Cor. to Prop. 2.2]{Drinfeld74}. Finally, $e(\phi)\in w(L^\times)$ 
	holds after a suitable ramified extension, which gives the last assertion. 
\end{proof}

By \cite[Prop. 7.1]{Drinfeld74}, $\phi$ has stable reduction after a finite extension; Lemma \ref{lemStable} makes this explicit.

\begin{lem}\label{lemIsotriv}
	\begin{enumerate}
		\item If $\phi\cong\phi'$ over some extension of $K$, then $J_i(\phi)=J_i(\phi')$ for all $i$.
		\item $\phi$ is isotrivial if and only if all $J_i(\phi)$ are algebraic over $\F_q$.
		\item $\phi$ has potentially good reduction at $v$ if and only if $v(J_i(\phi))\geq0$ for all $i$, if and only if $r_v(\phi)=r$.
		\item Assume that $K$ is finitely generated over $\F_q$ of transcendence degree $n\geq1$ and that $\phi$ is not isotrivial.
		Then there is a discrete valuation $v$ on $K$, trivial on $\F_q$, whose residue field is finitely generated over $\F_q$ of
		transcendence degree $n-1$, such that $\phi$ does not have potentially good reduction at $v$.
	\end{enumerate}
\end{lem}
\begin{proof}
	(1) Every isomorphism has the form $\phi'=c^{-1}\phi c$, which replaces $g_i$ by $c^{q^i-1}g_i$. This substitution does not change $J_i$.

	(2) If $\phi$ is isomorphic to a Drinfeld module $\psi$ over a finite field, then $J_i(\phi)=J_i(\psi)$ is algebraic over $\F_q$ by (1).
	Conversely, assume that all $J_i(\phi)$ are algebraic over $\F_q$. Choose $c\in\oK$ with $c^{q^N-1}=g_N^{-1}$, and put $\psi=c^{-1}\phi c$. 
	Then $\psi_a=\gamma(a)+\sum_ic^{q^i-1}g_i\tau^i$ has leading coefficient $1$, and $(c^{q^i-1}g_i)^{q^N-1}=J_i(\phi)$. So all coefficients 
	of $\psi_a$ lie in a finite field $k_0\subset\oK$. Let $b\in A$ and write $\psi_b=\sum_{j=0}^nu_j\tau^j$. Comparing the coefficients 
	of $\tau^{j+N}$ in $\psi_b\psi_a=\psi_a\psi_b$ gives 
	$$
	u_j^{q^N}-u_j=P_j(u_{j+1},\dots,u_n)
	$$ 
	for a polynomial $P_j$ with coefficients in $k_0$. By descending induction on $j$, every $u_j$ is algebraic over $k_0$. Applying this to 
	a finite set of generators $b$ of the $\F_q$-algebra $A$, we see that $\psi$ is defined over a finite field.

	(3) The conditions $v(J_i)\geq 0$ for all $i$ say that $v(g_N)/(q^N-1)\leq v(g_i)/(q^i-1)$ for all $i$, i.e., $N_v(\phi)=N$, i.e., 
	$r_v(\phi)=r$. This is equivalent to potentially good reduction by Lemma \ref{lemStable}.

	(4) By (2), some $J\colonequals J_i(\phi)$ is transcendental over $\F_q$. Extend it to a transcendence basis $J, x_2,\dots,x_n$ of 
	$K$ over $\F_q$, and put $E_0=\F_q(x_2,\dots,x_n)$, so that $K$ is a finite extension of $E=E_0(J)$. Let $v_0$ be the discrete 
	valuation of $E$ which is trivial on $E_0$ and satisfies $v_0(J)=-1$, i.e., $v_0(f/g)=\deg_J(g)-\deg_J(f)$ for $f,g\in E_0[J]$. 
	Its residue field is $E_0$. Let $v$ be any extension of $v_0$ to $K$. Then $v$ is a discrete valuation, trivial on $\F_q$, 
	whose residue field is a finite extension of $E_0$, hence finitely generated over $\F_q$ of transcendence degree $n-1$. 
	Since $v(J)<0$, the module $\phi$ does not have potentially good reduction at $v$ by (3). 
\end{proof}

\begin{rem}
	Let $K$ be finitely generated over $\F_q$. Then Lemma \ref{lemIsotriv} (2)--(4) imply that $\phi$ is isotrivial if and only 
	if it has potentially good reduction at every discrete valuation of $K$ trivial on $\F_q$. This is also a consequence of the 
	affineness of the coarse moduli space of Drinfeld modules; cf. the proof of \cite[Prop. 3.6]{PinkII}. For general $K$ the 
	``if'' part can fail, since $K$ may have no nontrivial discrete valuations.
\end{rem}

\begin{lem}\label{lemRedEnd}
	Assume that $\phi$ and $\phi'$ have coefficients in $R$, and that their reductions $\bar\phi$ and $\bar\phi'$ are Drinfeld modules
	over $k$. Then $\Hom_K(\phi,\phi')\subset\twist{R}$, and reduction modulo $\cM$ gives an injective homomorphism
	$$
	\Hom_K(\phi,\phi')\longhookrightarrow\Hom_k(\bar\phi,\bar\phi').
	$$
\end{lem}
\begin{proof}
	This is \cite[Lem. 3.3 (i), (ii)]{Gekeler83}, which is stated for an arbitrary $A$-field with a rank one valuation such that 
	$\gamma(A)\subset R$. For $A=\F_q[T]$ see also \cite[Lem. 6.1.9, Props. 6.1.10 and 6.1.11]{PapikianGTM}. 
\end{proof}


\section{Tate uniformization}\label{sTate}

Let $k$ be a finite $A$-field with $\chr_A(k)=\fp$. Let $K=\ls{k}{x}$, with $v=\ord_x$, $R=\dvr{k}{x}$, and $\abs{\cdot}=(\#k)^{-v(\cdot)}$,
viewed as an $A$-field via $A\to k\subset K$. Let $C$ be the completion of an algebraic closure of $K$. Let $\psi$ be a Drinfeld module 
of rank $r'$ over $k$. A \textit{$\psi$-lattice} of rank $m$ is a projective $A$-submodule $\La\subset{^\psi}K^\sep$ of rank $m$, stable under 
$\Gal(K^\sep/K)$, which has only finitely many elements in every disc of finite radius \cite[$\S$7]{Drinfeld74}. Let
$$
e_\La(x)=x\prodprime_{\la\in\La}\left(1-\frac{x}{\la}\right)
$$
be the associated exponential function. Its coefficients lie in $K$, and $e_\La$ is invertible in $\twistdvr{K}$ since its constant term is $1$.

Drinfeld's Tate uniformization \cite[Prop. 7.2]{Drinfeld74} is stated for an arbitrary complete discretely valued ring over $A$, 
so it applies in special characteristic. The description of endomorphisms is due to Gekeler \cite[Prop. 3.8]{Gekeler83}, 
again without any assumption on the characteristic. We record what we need. 

\begin{lem}\label{lemTate}
	\begin{enumerate}
		\item There is a unique Drinfeld module $\phi$ over $K$ such that $e_\La\psi_a=\phi_ae_\La$ in $\twistdvr{K}$ for all $a\in A$.
		It has rank $r'+m$, it has coefficients in $R$, and $\bar\phi=\psi$.
		\item The map $w\mapsto e_\La^{-1}we_\La$ is an isomorphism of $A$-algebras
		$$
		\End_C(\phi)\overset{\sim}{\To}\{u\in\End_C(\psi)\mid u\La\subseteq\La\}.
		$$
	\end{enumerate}
\end{lem}
\begin{proof}
	(1) This is the construction in \cite[Prop. 7.2]{Drinfeld74}; see also \cite[(3.6)]{Gekeler83}. Since every nonzero $\la\in\La$ 
	satisfies $\abs{\la}>1$, we have $e_\La\in1+\twistdvr{\cM}\tau$, which gives $\phi_b\in\twist{R}$ for all $b\in A$ and $\bar\phi=\psi$. 
	
	(2) This is \cite[Prop. 3.8]{Gekeler83}. In the proof there, the inclusion of the right-hand side into $\End_C(\phi)$ is said 
	to follow from the construction. For the reader's convenience we note how to see that $e_\La ue_\La^{-1}$ is a polynomial when 
	$\partial u=0$, which can happen only in special characteristic. If $\partial u\neq0$, the zeros of $e_\La(u(x))$ are simple and form 
	$u^{-1}\La\supseteq\La$, so $e_\La u=we_\La$ with $w\in\twist{C}$ by \cite[Lem. 5.1.4]{PapikianGTM}, whose proof works over any complete algebraically closed field of characteristic $p$. In general, choose 
	$c\in\{0,1\}\subset\F_q$ with $\partial(u+c)\neq0$. Then $u+c$ also preserves $\La$, and $e_\La ue_\La^{-1}=e_\La(u+c)e_\La^{-1}-c$. 
\end{proof}

\section{Proof of Theorem \ref{thmAp}}\label{sProof}

\subsection{From $\phi$ to $\psi$} Let $\phi$ be a non-isotrivial Drinfeld module of rank $r$ over an $A$-field $K$ with $\chr_A(K)=\fp$,
and let $D=\End^\circ(\phi)$. Since $A$ is a finitely generated $\F_q$-algebra, $\phi$ is defined over a subfield $K_0\subseteq K$ which is finitely generated over $\F_q$. Then $\End^\circ(\phi)$ is the same over $K$ and over $K_0$, and $\phi$ is
not isotrivial over $K_0$ by Lemma \ref{lemIsotriv} (2). So we may assume that $K$ is finitely generated over $\F_q$. We argue by induction
on $n=\mathrm{trdeg}_{\F_q}K$. If $n=0$, then $\phi$ is defined over a finite field, which is excluded. Let $n\geq1$.

Let $v$ be a discrete valuation of $K$ as in Lemma \ref{lemIsotriv} (4). Choose a finite separable extension $K'/K$ such that
$\End_{\oK}(\phi)=\End_{K'}(\phi)$, and extend $v$ to $K'$. Then there are a finite extension $L/K'$, an extension $w$ of $v$ to $L$, and
$c\in L^\times$ such that $\phi'=c^{-1}\phi c$ has coefficients in the valuation ring $R_w$ and its reduction $\bar\phi'$ is a Drinfeld module
over the residue field $k_w$ of rank $r_1=r_v(\phi)$. By Lemma \ref{lemIsotriv} (3), $r_1<r$. By Lemma \ref{lemRedEnd},
$$
\End_{\oK}(\phi)=\End_{K'}(\phi)\subseteq\End_L(\phi)\cong\End_L(\phi')\longhookrightarrow\End_{k_w}(\bar\phi'),
$$
so $D$ embeds into $D_1\colonequals\End^\circ(\bar\phi')$. Here $k_w$ is a finite extension of the residue field of $v$, so it is finitely
generated over $\F_q$ of transcendence degree $n-1$.

If $\bar\phi'$ is isotrivial, then it is isomorphic over $\overline{k_w}$ to a Drinfeld module $\psi$ over $\oF_\fp$, and
$D_1\cong\End^\circ(\psi)$. Put $r'=r_1$. Otherwise, the induction hypothesis applied to $\bar\phi'$ gives a Drinfeld module $\psi$ over
$\oF_\fp$ of rank $r'<r_1$ with $D_1\hookrightarrow\End^\circ(\psi)$. In both cases $D\hookrightarrow\End^\circ(\psi)$ with $r'<r$, and
$r\equiv r'\Mod{[D:F]}$ by Corollary \ref{corRankCong}. This proves (1) $\Rightarrow$ (2). The last assertion of Theorem \ref{thmAp} is
Corollary \ref{corRankCong}. That one can take $r'=r-[D:F]$ in (2) is shown in the proof of Theorem \ref{thmA} in $\S$\ref{sFF}. \qed

\subsection{From $\psi$ to $\phi$} Let $\psi$ be a Drinfeld module of rank $r'<r$ over $\oF_\fp$ with an embedding
$D\hookrightarrow\End^\circ(\psi)$, and let $n=(r-r')/[D:F]\geq1$. Choose a finite field $k\supseteq\F_\fp$ such that $\psi$ is defined over
$k$ and $\End_k(\psi)=\End_{\oF_\fp}(\psi)$. We identify $D$ with its image. Then
$$
\cO\colonequals D\cap\End_k(\psi)
$$
is an $A$-order in $D$, a projective $A$-module of rank $[D:F]$. Its elements are $\F_q$-linear polynomials with coefficients in $k$, so they act on
$K=\ls{k}{x}$. Let $\alpha_1,\dots,\alpha_n$ be distinct prime numbers different from $p$, put $\omega_i=x^{-\alpha_i}$, and let
$$
\La\colonequals\sum_{i=1}^n\cO(\omega_i)\subset K.
$$

\begin{lem}\label{lemNoCancel}
	If $u_1,\dots,u_n\in\twist{k}$ are not all zero, then
	$$
	\abs{\textstyle\sum_{i=1}^nu_i(\omega_i)}=\max_i\abs{u_i(\omega_i)}>1.
	$$
\end{lem}
\begin{proof}
	If $u\in\twist{k}$ is nonzero of degree $m$, then $\abs{u(\omega_i)}=\abs{\omega_i}^{q^m}=(\#k)^{\alpha_iq^m}$, because
	$\abs{\omega_i}>1$ and the nonzero coefficients of $u$ are units. For $i\neq j$ and any $m,m'\geq0$ we have $\alpha_iq^m\neq\alpha_jq^{m'}$
	by unique factorization. Hence the nonzero terms of the sum have pairwise distinct absolute values.
\end{proof}

By Lemma \ref{lemNoCancel}, the map $\cO^n\to\La$, $(u_i)\mapsto\sum_iu_i(\omega_i)$, is an isomorphism of $A$-modules, where $A$ acts
through $\psi(A)\subseteq\cO$. So $\La\cong\cO^n$ is a projective $A$-module of rank $n[D:F]$. It is discrete: if $\abs{\sum_iu_i(\omega_i)}\leq M$ with $u_i\in\cO$, then $\abs{u_i(\omega_i)}\leq M$ for all $i$ by 
Lemma \ref{lemNoCancel}, which bounds the degrees of the $u_i$, and $\twist{k}$ has only finitely many elements of bounded 
degree. It is Galois stable, because $\La\subset K$. Thus $\La$ is a $\psi$-lattice of rank $n[D:F]$. Let $\phi$ be the Drinfeld
module of rank $r'+n[D:F]=r$ over $K$ given by Lemma \ref{lemTate} (1).

Since $\bar\phi=\psi$ has rank $r'<r$, we have $r_v(\phi)=r'<r$ by Lemma \ref{lemStable}. By Lemma \ref{lemIsotriv} (3) and (2), $\phi$ is not isotrivial.

Since $\psi$ is defined over $k$, we have $\End_C(\psi)=\End_{\oF_\fp}(\psi)=\End_k(\psi)$ (see $\S$\ref{ssNotation}), 
and similarly $\End_C(\phi)=\End_{\oK}(\phi)$. Let $u\in\End_k(\psi)$ with $u\La\subseteq\La$. Then $u(\omega_1)=\sum_is_i(\omega_i)$ with $s_i\in\cO$, i.e., 
$(u-s_1)(\omega_1)-\sum_{i\geq2}s_i(\omega_i)=0$. By Lemma \ref{lemNoCancel}, $u=s_1\in\cO$. Conversely, $\cO$ preserves $\La$. 
By Lemma \ref{lemTate} (2), $\End_{\oK}(\phi)\cong\cO$, and therefore $\End^\circ(\phi)\cong D$.

Finally, $\phi$ is defined over a finitely generated subfield of $K$, and its endomorphism algebra over $\oK$ is the same as over the algebraic closure of this subfield. This proves (2) $\Rightarrow$ (1). \qed

\subsection{Endomorphism rings} The arguments above also classify the endomorphism rings. 

\begin{thm}\label{thmRings}
	Let $\fp\lhd A$ be a nonzero prime, let $\cO$ be an $A$-algebra, and let $r\geq1$. The following are equivalent: 
	\begin{enumerate}
		\item There are an $A$-field $K$ with $\chr_A(K)=\fp$ and a non-isotrivial Drinfeld module $\phi$ of rank $r$ over $K$ 
		such that $\End_{\oK}(\phi)\cong\cO$ as $A$-algebras. 
		\item There are a Drinfeld module $\psi$ of rank $r'<r$ over $\oF_\fp$ and an injective homomorphism of $A$-algebras 
		$\cO\hookrightarrow\End_{\oF_\fp}(\psi)$ such that $\dim_F(\cO\otimes_AF)$ divides $r-r'$. 
	\end{enumerate}
\end{thm}
\begin{proof}
	In (2), $\cO$ is a finitely generated torsion-free $A$-module, so $D\colonequals\cO\otimes_AF$ is a finite-dimensional $F$-subalgebra 
	of the division algebra $\End^\circ(\psi)$, hence a division $F$-algebra, and $\cO$ is an $A$-order in $D$. 
	
	(1) $\Rightarrow$ (2). The maps in the proof of (1) $\Rightarrow$ (2) of Theorem \ref{thmAp} are injective homomorphisms of $A$-algebras 
	$\End_{\oK}(\phi)\hookrightarrow\End_{\overline{k_w}}(\bar\phi')$, and in the isotrivial case $\End_{\overline{k_w}}(\bar\phi')\cong\End_{\oF_\fp}(\psi)$. 
	By induction we obtain $\cO\hookrightarrow\End_{\oF_\fp}(\psi)$ with $r'<r$, and the divisibility follows from Corollary \ref{corRankCong}. 
	
	(2) $\Rightarrow$ (1). Choose a finite field $k$ with $\psi$ defined over $k$ and $\End_k(\psi)=\End_{\oF_\fp}(\psi)$, and identify $\cO$ 
	with its image in $\End_k(\psi)$. The construction in the proof of (2) $\Rightarrow$ (1) of Theorem \ref{thmAp}, applied to this order $\cO$, 
	produces a non-isotrivial $\phi$ of rank $r$ with $\End_{\oK}(\phi)\cong\cO$. Indeed, the argument there used only that $\cO$ 
	is an $A$-order in $D$ contained in $\End_k(\psi)$. 
\end{proof}

\begin{rem}
	\begin{enumerate}
		\item The case $\cO=\End_k(\psi)$ of the construction gives, for every Drinfeld module $\psi$ of rank $r'$ over a finite field
		and every $n\geq1$, a non-isotrivial Drinfeld module of rank $r'+n\dim_F\End^\circ(\psi)$ whose endomorphism ring over $\oK$ is 
		isomorphic to $\End_{\oF_\fp}(\psi)$.
		\item Theorem \ref{thmAp} also controls the endomorphism algebras over $K$ itself. If $\phi$ is non-isotrivial of rank $r$, then
		$\End^\circ_K(\phi)\subseteq\End^\circ(\phi)$ embeds into $\End^\circ(\psi)$ for some $\psi$ over $\oF_\fp$ of rank $r'<r$, and
		$r\equiv r'\Mod{[\End^\circ_K(\phi):F]}$ by Corollary \ref{corRankCong}.
	\end{enumerate}
\end{rem}


\section{Subalgebras over finite fields and proof of Theorem \ref{thmA}}\label{sFF}

Let $\psi$ be a Drinfeld module of rank $r'$ over a finite field $k$ with $\chr_A(k)=\fp$, and let $\pi$ be its Frobenius. Put 
$$
\Delta\colonequals\End^\circ_k(\psi),\qquad E\colonequals F(\pi),\qquad h\colonequals[E:F],\qquad m\colonequals r'/h. 
$$
By \cite[Thm. 2.9]{GekelerFDM}, \cite[Thm. 1]{Yu95}, \cite[(2.2.2)]{Laumon}, $\Delta$ is a central division algebra over $E$ of 
degree $m$, the field $E$ has a unique place $\infty_E$ over $\infty$, there is a unique place $\fP$ of $E$ dividing $\pi$, it lies 
over $\fp$, and 
$$
\inv_{\infty_E}(\Delta)=-1/m,\qquad \inv_\fP(\Delta)=1/m,\qquad \inv_v(\Delta)=0\quad\text{for } v\neq\infty_E,\fP. 
$$
Every $\psi$ over $\oF_\fp$ is defined over such a $k$ with $\End^\circ(\psi)=\End^\circ_k(\psi)$. 

\begin{prop}\label{propFF}
	Let $D$ be a division $F$-algebra with center $Z$, and put $e=e(D)$, $d=d(D)$. Then $D$ embeds into $\Delta$ as an $F$-algebra 
	if and only if there is a field $L$ which is a quotient of $Z\otimes_FE$ with the following properties: 
	\begin{enumerate}
		\item $L$ has a unique place $\infty_L$ over $\infty_E$ and a unique place $\fP_L$ over $\fP$, and $[L:F]\,d$ divides $r'$; 
		\item $D_L\colonequals D\otimes_ZL$ splits at every place of $L$ other than $\infty_L$ and $\fP_L$, and 
		$s_L\cdot\inv_{\infty_L}(D_L)=-1/d$ in $\Q/\Z$, where $s_L\colonequals r'/([L:F]d)$.
	\end{enumerate}
\end{prop}
\begin{proof}
	For a field $L\supseteq E$ put $\ell\colonequals[L:E]$. By (C6), $L$ is $E$-isomorphic to a subfield of $\Delta$ 
	if and only if $\ell$ divides $m$ and $\frac{m}{\ell}[L_w:E_v]\inv_v(\Delta)=0$ for all places $w\mid v$ of $L$ and $E$. Since 
	$[L_w:E_v]\leq\ell$, this holds if and only if 
	\begin{equation}\label{eqLinDelta}
		\ell\mid m,\quad\text{and $L$ has a unique place $\infty_L$ over $\infty_E$ and a unique place $\fP_L$ over $\fP$.}
	\end{equation}
	If $L\subseteq\Delta$, then $\Delta_L\colonequals\Cent_\Delta(L)$ is a central division algebra over $L$ of degree 
	$n\colonequals m/\ell=r'/[L:F]$, Brauer equivalent to $\Delta\otimes_EL$ by (C2). So by (C5) its invariants are $-1/n$ at $\infty_L$, $1/n$ at 
	$\fP_L$, and $0$ elsewhere. 

	Suppose $D\subseteq\Delta$. The subring $L\colonequals ZE$ of $\Delta$ is commutative and finite-dimensional over $F$, hence a field, 
	and it is a quotient of $Z\otimes_FE$. It satisfies \eqref{eqLinDelta}. The algebra $D_L=D\otimes_ZL$ is central simple 
	over $L$, so the natural map $D_L\to\Delta$ is injective, and its image lies in $\Delta_L$. Thus $D_L$ is a division 
	$F$-subalgebra of $\End^\circ(\psi)$ with center $L$. By Proposition \ref{propInvInf}, $[L:F]d$ divides $r'$, 
	$s_L\cdot\inv_{\infty_L}(D_L)=-1/d$, and $\gcd(d,s_L)=1$. Moreover 
	$\Delta_L\cong D_L\otimes_L\Cent_{\Delta_L}(D_L)$ by (C1), where the second factor has degree $n/d=s_L$. At a place $w\neq\infty_L,\fP_L$ of 
	$L$ we have $\inv_w(\Delta_L)=0$, so $\inv_w(D_L)$ is annihilated by both $d$ and $s_L$, hence is $0$. This proves (1) and (2). 

	Conversely, let $L$ satisfy (1) and (2). Then $m/\ell=r'/[L:F]=ds_L$ is an integer, so \eqref{eqLinDelta} holds and we may 
	identify $L$ with a subfield of $\Delta$. Put $\iota\colonequals\inv_{\infty_L}(D_L)$. By (2) and (C4), $D_L$ has invariants 
	$\iota$ at $\infty_L$, $-\iota$ at $\fP_L$, and $0$ elsewhere. Hence $\Delta_L\otimes_LD_L^{\mathrm{opp}}$ has invariants 
	$$
	\mp\Bigl(\frac{1}{ds_L}+\iota\Bigr)\quad\text{at }\infty_L,\ \fP_L, 
	$$
	and $0$ elsewhere. Since $s_L\bigl(\frac{1}{ds_L}+\iota\bigr)=\frac1d-\frac1d=0$, all of them lie in $\frac{1}{s_L}\Z/\Z$, so by (C4) 
	the algebra $\Delta_L\otimes_LD_L^{\mathrm{opp}}$ is 
	Brauer equivalent to a central division algebra $B$ over $L$ whose degree $b$ divides $s_L$. The algebras $\Delta_L$ and 
	$D_L\otimes_L\Mat_{s_L/b}(B)$ have the same degree $ds_L$ and the same Brauer class, so they are isomorphic. Therefore 
	$D\subseteq D_L\hookrightarrow\Delta_L\subseteq\Delta$. 
\end{proof}

Proposition \ref{propFF} depends on $\psi$ through $E$ and $\fP$. Letting $\psi$ vary, we obtain a criterion in terms of $D$ alone. 

\begin{cor}\label{corFFinv}
	Let $D$ be a division $F$-algebra with center $Z$, put $e=e(D)$, $d=d(D)$, and let $r'\geq1$. There is a Drinfeld module $\psi$ of 
	rank $r'$ over $\oF_\fp$ with an embedding of $F$-algebras $D\hookrightarrow\End^\circ(\psi)$ if and only if the following hold: 
	\begin{enumerate}
		\item[(i)] $ed$ divides $r'$; 
		\item[(ii)] $Z$ has a unique place $\infty_Z$ over $\infty$; 
		\item[(iii)] there is a place $\fp_Z$ of $Z$ over $\fp$ such that $D$ splits at every place of $Z$ other than $\infty_Z$ and $\fp_Z$; 
		\item[(iv)] $s'\cdot\inv_{\infty_Z}(D)=-1/d$ in $\Q/\Z$, where $s'=r'/(ed)$. 
	\end{enumerate}
\end{cor}
\begin{proof}
	Suppose $D\subseteq\End^\circ(\psi)$. Then (i), (ii) and (iv) hold by Proposition \ref{propInvInf}, which also gives $\gcd(d,s')=1$. 
	Choose $k$ as above with $\End^\circ(\psi)=\End^\circ_k(\psi)$, let $L$ be as in Proposition \ref{propFF}, and put 
	$\fp_Z=\fP_L\cap Z$, $f=[L:Z]$. Then $s'=fs_L$, so $\gcd(d,f)=1$. Let $u$ be a place of $Z$ other than $\infty_Z$ and $\fp_Z$. 
	No place of $L$ over $u$ equals $\infty_L$ or $\fP_L$, so for every $w\mid u$ we have 
	$0=\inv_w(D_L)=[L_w:Z_u]\cdot\inv_u(D)$ by Proposition \ref{propFF} (2) and (C5). Summing over $w\mid u$ and using 
	$\sum_{w\mid u}[L_w:Z_u]=f$, we get $f\cdot\inv_u(D)=0$. Since also $d\cdot\inv_u(D)=0$ and $\gcd(d,f)=1$, we conclude that 
	$\inv_u(D)=0$. This proves (iii). 
	
	Conversely, assume (i)--(iv). Let $\cO_Z$ be the integral closure of $A$ in $Z$. Its class group is finite, so there are $j\geq1$ 
	and $\varpi\in\cO_Z$ with $\varpi\cO_Z=\fp_Z^j$. Let $E=F(\varpi)\subseteq Z$, $h=[E:F]$, and $\fP=\fp_Z\cap E$. Then $\varpi$ is 
	integral over $A$, and $\infty$ has a unique place in $E$ by (ii). The only zero of $\varpi$ in $Z$ is $\fp_Z$, because 
	$\varpi\cO_Z=\fp_Z^j$ and $\varpi$ has a pole at $\infty_Z$. Hence $\fP$ is the only zero of $\varpi$ in $E$, and $\fp_Z$ is the 
	only place of $Z$ over $\fP$, since $\varpi$ vanishes at every such place. Put $N_0=\ord_\fP(\varpi)\,[\F_\fP:\F_\fp]$, and let $k\supseteq\F_\fp$ be the finite field with 
	$[k:\F_\fp]=N_0r'/h$, an integer since $h\mid e\mid r'$. Since $\infty$ does not split in $E$ and $\Nr_{E/F}(\varpi)$ 
	generates $\fp^{N_0}$, the extension of $\abs{\cdot}$ to $E\otimes_F\Fi$ satisfies 
	$$
	\abs{\varpi}^h=\abs{\Nr_{E/F}(\varpi)}=(\#\F_\fp)^{N_0}=(\#k)^{h/r'}. 
	$$
	Together with $h\mid r'$, this shows that $\varpi$ satisfies conditions (1)--(5) of \cite[p. 165]{Yu95}, i.e., it is a Weil number of rank $r'$ over $k$. By \cite[Thm. 3]{Yu95} (see also 
	\cite[(2.2.2) (iv)]{Laumon}), there is a Drinfeld module $\psi$ of rank $r'$ over $k$ whose Frobenius is conjugate to $\varpi$ 
	over $F$. We identify $F(\pi)$ with $E$ so that $\pi=\varpi$. Then $\fP$ is the place of $E$ dividing $\pi$. 
	Proposition \ref{propFF} applies to $\psi$ and $L=Z$, with $D_L=D$ and $s_L=s'$: condition (1) holds by (i) and the choice of 
	$\varpi$, and condition (2) is (iii) and (iv). Hence $D\hookrightarrow\End^\circ_k(\psi)\subseteq\End^\circ(\psi)$. 
\end{proof}

\begin{rem}\label{remFFsuper}
	Conditions (ii)--(iv) do not require $\fp$ to have a unique place in $Z$. In Proposition \ref{propFF} only $\fP$ is required 
	not to split in $L$, while other places of $E$ over $\fp$ may split; see Remark \ref{remExOtherChar} for an example with $L=Z=E$. 
	If $\fp$ does have a unique place in $Z$, then in the proof of Corollary \ref{corFFinv} one can take $\varpi\in A$ and $\psi$ 
	supersingular. 
\end{rem}

\begin{proof}[Proof of Theorem \ref{thmA}]
	By Theorem \ref{thmAp} and Corollary \ref{corFFinv}, there is a non-isotrivial $\phi$ of rank $r$ and characteristic $\fp$ with 
	$\End^\circ(\phi)\cong D$ if and only if conditions (i)--(iv) of Corollary \ref{corFFinv} hold for some integer $r'$ with $0<r'<r$ and 
	$r'\equiv r\Mod{ed^2}$. 

	Suppose that such an $r'$ exists, and put $s'=r'/(ed)$. Since $ed$ divides $r'$ and $r-r'$, it divides $r$. Moreover, 
	$s'\equiv s\Mod{d}$ and $d\cdot\inv_{\infty_Z}(D)=0$, so (iv) is equivalent to (4), and $0<s'<s$ implies $s\geq s'+d>d$, i.e., $[D:F]=ed^2<eds=r$. 
	Thus (1)--(4) hold. Conversely, if (1)--(4) hold, then $r'\colonequals r-ed^2=ed(s-d)$ satisfies $0<r'<r$ and 
	$r'/(ed)\equiv s\Mod{d}$, so (i)--(iv) hold for $r'$. The last assertion of Theorem \ref{thmA} is part of Theorem \ref{thmAp}. 
\end{proof}

The proof shows that if (1) of Theorem \ref{thmAp} holds, then one can take $r'=r-[D:F]$ in (2), as asserted there. 

\begin{cor}\label{corSplit}
	Let $\phi$ be a Drinfeld module over an $A$-field of characteristic $\fp\neq0$, and let $D\subseteq\End^\circ(\phi)$ be a division 
	$F$-subalgebra with center $Z$. Then there is a place $\fp_Z$ of $Z$ over $\fp$ such that $D$ splits at every place of $Z$ other 
	than $\fp_Z$ and the place over $\infty$. This is Theorem \ref{thmB} (2) in special characteristic. 
\end{cor}
\begin{proof}
	If $\phi$ is isotrivial, then $\End^\circ(\phi)\cong\End^\circ(\psi)$ for a Drinfeld module $\psi$ over $\oF_\fp$. Otherwise 
	$\End^\circ(\phi)$ embeds into such an algebra by Theorem \ref{thmAp}. In both cases Corollary \ref{corFFinv} (iii) applies to $D$. 
\end{proof}

\begin{cor}\label{corPink}
	Let $\phi$ be a non-isotrivial Drinfeld module of rank $r$ over an $A$-field of characteristic $\fp\neq0$, let $D=\End^\circ(\phi)$, 
	and put $e=e(D)$, $d=d(D)$. 
	\begin{enumerate}
		\item $ed^2<r$. In particular $ed<r$, i.e., every commutative subfield of $D$ has degree $<r$ over $F$. 
		\item If $D$ is not commutative, then $r\geq ed(d+1)\geq6$. 
	\end{enumerate}
	In particular, the endomorphism algebra of a non-isotrivial Drinfeld module of rank $\leq5$ is a field. 
\end{cor}
\begin{proof}
	(1) is Theorem \ref{thmA} (1). In (2) we have $d\geq2$, and $s=r/(ed)$ is an integer with $ed^2<eds$, so $s\geq d+1$. 
\end{proof}

\begin{rem}\label{remPinkKummer}
	\begin{enumerate}
		\item The inequality $e(D)d(D)<r$ is \cite[Thm. 6.1]{PinkII}, where it is deduced from the finiteness of the moduli stack of 
		rank one Drinfeld modules. Since $e(D)d(D)=r$ when $\phi$ is isotrivial, Corollary \ref{corPink} (1) shows that $e(D)d(D)=r$ if and 
		only if $\phi$ is isotrivial, which is the equivalence (a) $\Leftrightarrow$ (e) in \cite[Prop. 6.3]{DevicPink}. 
		\item Suppose $R=\End_{\oK}(\phi)$ has center $A$, so $e=1$, and write $r=mn$ with $m=d$, $n=s$ as in 
		\cite[$\S$4]{PinkKummer}. For a maximal order $R$, Pink writes $R\otimes_AA_\fl\cong\Mat_{n_\fl}(S_\fl)$ with $S_\fl$ a maximal 
		order in a central division algebra over $F_\fl$ that may depend on $\fl$ \cite[$\S$5]{PinkKummer}. Corollary \ref{corSplit} 
		shows that $S_\fl=A_\fl$ for every prime $\fl\neq\fp$, and Proposition \ref{propInvInf} shows that $\gcd(m,n)=1$. 
	\end{enumerate}
\end{rem}

For a Drinfeld module $\phi$ of rank $r$ over an $A$-field, we call 
$$
\mu(\phi)\colonequals\frac{r}{[\End^\circ(\phi):F]}
$$
the \textit{endomorphism ratio} of $\phi$. It is the analogue of the ratio $2\dim(X)/[\End^0(X):\Q]$ of a simple abelian 
variety $X$ considered by Oort \cite[Question 22.2]{OortQ}. 

\begin{cor}\label{corRatio}
	Let $\fp\lhd A$ be a nonzero prime. 
	\begin{enumerate}
		\item If $\phi$ is isotrivial of $A$-characteristic $\fp$, then $\mu(\phi)=1/d$, where $d=d(\End^\circ(\phi))$. 
		Every $d\geq1$ occurs. 
		\item If $\phi$ is non-isotrivial of $A$-characteristic $\fp$ and $D=\End^\circ(\phi)$, then $\mu(\phi)=s/d$ with 
		$d=d(D)$, $s=r/(e(D)d)$, $\gcd(s,d)=1$ and $s>d$. In particular $\mu(\phi)>1$. 
		\item Let $\mu>1$ be a rational number, written as $\mu=s/d$ in lowest terms, and let $r\geq1$. There is a non-isotrivial 
		Drinfeld module of rank $r$ and $A$-characteristic $\fp$ with endomorphism ratio $\mu$ if and only if $ds$ divides $r$. 
	\end{enumerate}
	Thus the endomorphism ratios of Drinfeld modules of $A$-characteristic $\fp$ are the numbers $1/d$, $d\geq1$, and the 
	rational numbers $>1$, and $\mu(\phi)\leq1$ if and only if $\phi$ is isotrivial. 
\end{cor}
\begin{proof}
	(1) Here $\End^\circ(\phi)\cong\End^\circ_k(\psi)$ for a Drinfeld module $\psi$ of rank $r$ over a finite field $k$, and this is 
	a central division algebra over $F(\pi)$ of degree $r/[F(\pi):F]$. So $r=ed$ and $\mu(\phi)=r/(ed^2)=1/d$, where 
	$e=e(\End^\circ(\phi))$. Given $d\geq1$, let $D$ be the central division algebra over $F$ with invariants $-1/d$ at $\infty$, 
	$1/d$ at $\fp$, and $0$ elsewhere. By Corollary \ref{corFFinv}, $D$ embeds into $\End^\circ(\psi)$ for some $\psi$ of rank $d$ 
	over $\oF_\fp$. Since $[\End^\circ(\psi):F]\leq d^2$, we get $\End^\circ(\psi)=D$ and $\mu(\psi)=1/d$. 
	
	(2) We have $\mu(\phi)=r/(ed^2)=s/d$. Here $\gcd(s,d)=1$ by Theorem \ref{thmB} (1), and $s>d$ by Theorem \ref{thmA} (1). 
	
	(3) If such a $\phi$ exists, then $r=e(D)ds$ by (2), so $ds$ divides $r$. Conversely, let $r=eds$ with $e\geq1$. Let $u\in F$ 
	be a uniformizer at $\infty$, and let $Z=F(u^{1/e})$. The polynomial $x^e-u$ is Eisenstein at $\infty$, so $[Z:F]=e$ and $\infty$ 
	is totally ramified in $Z$. Choose a place $\fp_Z$ of $Z$ over $\fp$. Since $\gcd(s,d)=1$, there is a unique 
	$\iota\in\frac1d\Z/\Z$ with $s\iota=-1/d$, and it has order $d$. Let $D$ be the central division algebra over $Z$ with 
	invariants $\iota$ at $\infty_Z$, $-\iota$ at $\fp_Z$, and $0$ elsewhere; see (C4). Then $d(D)=d$, and $[D:F]=ed^2<eds=r$ because 
	$s>d$. By Theorem \ref{thmA}, $D\cong\End^\circ(\phi)$ for a non-isotrivial $\phi$ of rank $r$ and $A$-characteristic $\fp$, 
	and $\mu(\phi)=s/d$. 
\end{proof}

\begin{rem}\label{remRatio}
	\begin{enumerate}
		\item By Corollary \ref{corRatio} (3) with $\mu=(d+1)/d$, the bound $r\geq ed(d+1)$ of Corollary \ref{corPink} (2) is 
		attained for all $e\geq1$ and $d\geq2$. For $A=\F_q[T]$, $\fp=(T)$ and $e=1$, see Remark \ref{remExGeneral} with $j=1$. 
		In contrast, over finite fields noncommutative endomorphism algebras occur in every rank $\geq2$, and in generic 
		characteristic they do not occur at all. 
		\item In generic characteristic $\mu(\phi)=r/e$ is a positive integer, and every positive integer occurs. For abelian 
		varieties in characteristic $p$, every positive rational number is the endomorphism ratio of a simple abelian variety by 
		\cite{Tafazolian}, as expected by Oort \cite[p. 315]{OortQ}. The inequality $\mu(\phi)>1$ for non-isotrivial Drinfeld modules comes from the existence of a place of 
		bad reduction (Lemma \ref{lemIsotriv} (4)), which has no analogue for abelian varieties. 
	\end{enumerate}
\end{rem}


\section{Examples}\label{sExamples}

In this section $A=\F_q[T]$, $F=\F_q(T)$, $\infty$ is the place at which $1/T$ is a uniformizer, and $\fp=(T)$, so $\F_\fp=\F_q$. 
For an $A$-field $K$ we write $t=\gamma(T)$, and we use $a=T$ in the definition of the invariants $J_i$.

\begin{example}\label{ex1}
	Let $q$ be odd, let $c\in\F_q^\times$ be a non-square, and let $\beta\in\F_{q^2}$ with $\beta^2=c$. Let $B=\F_q[y]$, and embed $A$ into $B$
	by $T\mapsto y^2-c$. Then $Z\colonequals\F_q(y)=F(\sqrt{T+c})$ is a separable quadratic extension of $F$ with constant field $\F_q$.
	The prime $(T)=(y^2-c)$ is inert in $Z$, and $\infty$ is ramified in $Z$.

	Let $K=\F_{q^2}(x)$ with $x$ transcendental, viewed as an $A$-field via $\gamma(T)=0$. Let $\phi'\colon B\to\twist{K}$ be the Drinfeld
	$B$-module of rank $2$ defined by $\phi'_y=\beta+x\tau+\tau^2$, and let $\phi=\phi'|_A$. Using $\tau\beta=\beta^q\tau=-\beta\tau$ one computes
	$$
	\phi_T=(\phi'_y)^2-c=(2\beta+x^{q+1})\tau^2+(x+x^{q^2})\tau^3+\tau^4.
	$$
	So $\phi$ is a Drinfeld module of rank $4$ over $K$ of $A$-characteristic $(T)$. The element $u=\phi'_y$ lies in $\End_K(\phi)$ and
	satisfies $u^2=\phi_T+c$, so $Z\cong F(u)\subseteq\End^\circ(\phi)$.

	\begin{enumerate}
		\item \textit{$\phi$ is not isotrivial.} For the valuation $v$ of $K$ with $v(x)=-1$ we have $v(J_2(\phi))=-(q+1)(q^4-1)<0$.
		So $\phi$ does not have potentially good reduction at $v$ by Lemma \ref{lemIsotriv} (3), and $\phi$ is not isotrivial.

		\item \textit{$\End^\circ(\phi)=\End^\circ_K(\phi)=Z$, and $\End_{\oK}(\phi)=\F_q[u]\cong B$.} Let $D=\End^\circ(\phi)$. Since $Z\subseteq D$ is a 
		commutative subfield, $e(D)d(D)\geq 2$. By Theorem \ref{thmA} (1), $e(D)d(D)^2<4$. Hence $d(D)=1$ and $e(D)=2$, so $D$ is a 
		quadratic field containing $Z$, i.e., $D=Z$. The ring $\End_{\oK}(\phi)$ is an $A$-order in $Z$ containing 
		$\F_q[u]\cong B=\F_q[y]$, which is the maximal order of $Z$. So $\End_{\oK}(\phi)=\F_q[u]$. 
		
		\item \textit{$Z$ is not isomorphic to $\End^\circ_{k'}(\psi)$ for any Drinfeld module $\psi$ over a finite $A$-field $k'$ of 
		$A$-characteristic $(T)$.} 
		Suppose it is. Since $Z$ is commutative, $Z=F(\pi)$ by \cite[Thm. 4.1.3]{PapikianGTM}, where $\pi$ is the Frobenius of $k'$. 
		Moreover, $\pi$ has nonzero valuation only at the place of $Z$ over $\infty$ and at a place over $(T)$. These are the unique places 
		of $Z$ over $\infty$ and over $(T)$, so they are fixed by the nontrivial automorphism $\sigma$ of $Z/F$. Hence $\pi^\sigma/\pi$ 
		has trivial divisor, so $\pi^\sigma=\zeta\pi$ with $\zeta\in\F_q^\times$, and applying $\sigma$ again gives $\zeta^2=1$. 
		If $\zeta=1$, then $\pi\in F$, contradicting $F(\pi)=Z$. If $\zeta=-1$, then $\pi^2\in F$ has divisor supported at $(T)$ and $\infty$, 
		so $\pi^2=c'T^m$ with $c'\in\F_q^\times$ and $Z=F(\sqrt{c'T^m})$. If $m$ is odd, then $(T)$ is ramified in $Z$. If $m$ is even, 
		then $Z$ is a constant field extension. Both are impossible. 
	\end{enumerate}

	The field $Z$ satisfies the conditions of Theorem \ref{thmA} for $r=4$. In terms of Theorem \ref{thmAp}: since $(T)$ and $\infty$ do not split in $Z$, the field $Z$ embeds into the quaternion
	algebra over $F$ ramified at $(T)$ and $\infty$, which is the endomorphism algebra of a supersingular Drinfeld module of rank $2$
	\cite[Rem. 4.4.2 (2)]{PapikianGTM}. Here $r'=2$ and $[Z:F]=2$ divides $4-2$.
\end{example}

\begin{rem}\label{remExOtherChar}
	The restriction to characteristic $(T)$ in Example \ref{ex1} (3) cannot be dropped. Let $\fq\lhd A$ be a prime which splits in $Z$, 
	let $\fQ$ be a prime of $B$ above $\fq$, and let $\rho\colon B\to\twist{H}$ be a Drinfeld $B$-module of rank $1$ over a finite extension 
	$H$ of $Z$ with good reduction at a prime above $\fQ$ (for instance a sign-normalized module over the Hilbert class field, 
	cf. \cite[Ch. 7]{PapikianGTM}). Its restriction $\rho|_A$ is a Drinfeld $A$-module of rank $2$ with $\End^\circ(\rho|_A)\supseteq Z$. 
	The reduction $\bar\rho$ at a prime above $\fQ$ is a rank $1$ Drinfeld $B$-module over a finite field $k'$ of $B$-characteristic $\fQ$, 
	and $\bar\rho|_A$ has $A$-characteristic $\fq$ and height $[B_\fQ:A_\fq]\cdot H(\bar\rho)=1$. So $\bar\rho|_A$ is ordinary, 
	and $\End^\circ_{\overline{k'}}(\bar\rho|_A)$ is a quadratic field containing $Z$ \cite[Cor. 4.1.12]{PapikianGTM}, i.e., it equals $Z$. 
\end{rem}

\begin{example}\label{ex2}
	For $n\geq1$, let $D_n$ be the central division algebra over $F$ with $\inv_{(T)}(D_n)=1/n$ and $\inv_\infty(D_n)=-1/n$ which
	is split at all other places. By \cite[Rem. 4.4.2 (2)]{PapikianGTM}, $D_n\cong\End^\circ(\psi_n)$, where $(\psi_n)_T=\tau^n$. Let $D=D_3^{\mathrm{opp}}$, so that $\inv_{(T)}(D)=-1/3$ and
	$\inv_\infty(D)=1/3$.

	Let $K=\F_{q^3}(x)$ with $x$ transcendental, viewed as an $A$-field via $\gamma(T)=0$. Let $u=x\tau^2+\tau^5\in\twist{K}$, and 
	let $\phi$ be the Drinfeld module over $K$ defined by 
	$$
	\phi_T=u^3=x^{1+q^2+q^4}\tau^6+\bigl(x^{1+q^2}+x^{1+q^7}+x^{q^5+q^7}\bigr)\tau^9+\bigl(x+x^{q^5}+x^{q^{10}}\bigr)\tau^{12}+\tau^{15}. 
	$$
	It has rank $15$ and $A$-characteristic $(T)$; cf. \cite[Ex. 5.2 (c)]{PinkII} and \cite[Ex. 4.1.6]{PapikianGTM} for a similar example of rank $6$. 
	\begin{enumerate}
		\item \textit{$D$ is not isomorphic to $\End^\circ_{k'}(\psi)$ for any Drinfeld module $\psi$ over any finite $A$-field $k'$, of any $A$-characteristic.}
		The center of $D$ is $F$, so such a $\psi$ would have $F(\pi)=F$. Then $\rank(\psi)=3$, and $\End^\circ_{k'}(\psi)$ has invariant $-1/3$ at $\infty$
		by \cite[Thm. 4.1.3, Rem. 4.1.5 (2)]{PapikianGTM}. But $\inv_\infty(D)=1/3\neq-1/3$ in $\Q/\Z$.
		\item \textit{$\phi$ is not isotrivial.} For the valuation $v$ of $K$ with $v(x)=-1$ we have 
		$v(J_{12}(\phi))=-q^{10}(q^{15}-1)<0$. So $\phi$ does not have potentially good reduction at $v$ by 
		Lemma \ref{lemIsotriv} (3), and $\phi$ is not isotrivial. 
		\item \textit{$\End^\circ(\phi)=\End^\circ_K(\phi)\cong D$, and $\End_{\oK}(\phi)=\F_{q^3}[u]$.} For $\alpha\in\F_{q^3}$ we have 
		$u\alpha=\alpha^{q^2}u$, because $\alpha^{q^5}=\alpha^{q^2}$. Hence $\alpha$ commutes with $\phi_T=u^3$, and 
		$\cO\colonequals\F_{q^3}[u]$ is a subring of $\End_K(\phi)$. The map $\alpha\mapsto\alpha$, $u\mapsto\tau$ is an 
		anti-isomorphism of $\cO$ onto $\twist{\F_{q^3}}=\End_{\oF_q}(\psi_3)$ which maps $\phi_T$ to $(\psi_3)_T=\tau^3$. So 
		$\cO\otimes_AF\cong D_3^{\mathrm{opp}}=D$, and $\cO$ is a maximal $A$-order in this algebra by \cite[Thm. 4.3]{GekelerFDM}. 
		Thus $D$ embeds into $D'\colonequals\End^\circ(\phi)$, and $[D':F]$ is a multiple of $[D:F]=9$. Since $[D':F]<15$ by 
		Theorem \ref{thmA} (1), we get $D'=D$, and then $\End_{\oK}(\phi)=\cO$ by maximality. 
		\item \textit{Degeneration.} Let $v$ be as in (2), and let $c\in\oK$ with $c^{q^2-1}=x^{-1}$. Then 
		$c^{-1}uc=\tau^2+c^{q^5-1}\tau^5$ reduces to $\tau^2$, so $c^{-1}\phi c$ reduces to $\psi_6$, and the reduction map of 
		Lemma \ref{lemRedEnd} sends $\cO$ onto $\F_{q^3}[\tau^2]\subset\twist{\F_{q^6}}=\End_{\oF_q}(\psi_6)$. This is an embedding 
		$D\hookrightarrow D_6$ as in Theorem \ref{thmAp}, with $r'=6=15-[D:F]$. As a check, Proposition \ref{propInvInf} gives 
		$s=5$, and indeed $5\cdot\frac13=-\frac13$ in $\Q/\Z$. 
	\end{enumerate}
	More generally, Theorem \ref{thmA} with $e=1$, $d=3$ and $\fp_Z=(T)$ shows that $D_3$, for which $\inv_\infty=-1/3$, is the 
	endomorphism algebra of a non-isotrivial Drinfeld module of rank $r$ if and only if $r\equiv3\Mod9$ and $r\geq12$. Similarly, 
	$D_3^{\mathrm{opp}}$, for which $\inv_\infty=1/3$, occurs if and only if $r\equiv6\Mod9$ and $r\geq15$. So the rank of $\phi$ 
	is the smallest possible; see Remark \ref{remExGeneral}. 
\end{example}

\begin{rem}\label{remExGeneral}
	The construction of Example \ref{ex2} works for every central division algebra over $F$ which is ramified exactly at $(T)$ 
	and $\infty$. Let $d\geq2$ and $1\leq j<d$ with $\gcd(j,d)=1$. Let $K=\F_{q^d}(x)$, viewed as an $A$-field via $\gamma(T)=0$, 
	let $u=x\tau^j+\tau^{j+d}\in\twist{K}$, and let $\phi$ be the Drinfeld module of rank $r=d(j+d)$ over $K$ defined by 
	$\phi_T=u^d$. Then $\phi$ is not isotrivial, $\End_{\oK}(\phi)=\F_{q^d}[u]$, and $\End^\circ(\phi)$ is the central division 
	algebra $D$ over $F$ of dimension $d^2$ which is split away from $(T)$ and $\infty$ and satisfies 
	$j\cdot\inv_\infty(D)=-1/d$ in $\Q/\Z$. 

	Indeed, the coefficient of $\tau^{r-d}$ in $\phi_T$ is $\sum_{i=0}^{d-1}x^{q^{i(j+d)}}$, so $v(J_{r-d}(\phi))<0$ for the 
	valuation $v$ with $v(x)=-1$, and $\phi$ is not isotrivial by Lemma \ref{lemIsotriv} (3). For $\alpha\in\F_{q^d}$ we have 
	$u\alpha=\alpha^{q^j}u$, so $\cO\colonequals\F_{q^d}[u]$ is a subring of $\End_K(\phi)$. It is a twisted polynomial ring in 
	$u$ over $\F_{q^d}$, free of rank $d^2$ over $\phi(A)=\F_q[u^d]$, and its center is $\F_q[u^d]$ because 
	$\alpha\mapsto\alpha^{q^j}$ generates $\Gal(\F_{q^d}/\F_q)$. Hence $D\colonequals\cO\otimes_AF$ is a division subalgebra of 
	$D'\colonequals\End^\circ(\phi)$ with center $F$ and $[D:F]=d^2$. Now $[D':F]$ is a multiple of $d^2$, and 
	$[D':F]<r<2d^2$ by Theorem \ref{thmA} (1), so $D'=D$. The invariants of $D$ are given by Theorem \ref{thmB} (1) with 
	$s=j+d$, by Corollary \ref{corSplit}, and by the sum formula. Finally, let $w\in\End_{\oK}(\phi)$, and choose a nonzero 
	$f\in A$ with $\phi_fw\in\cO$. Since $\phi_f$ is central in $\cO$ and its leading coefficient is in $\F_q^\times$, division 
	with remainder in $\cO$ gives $\phi_fw=\phi_fh+\rho$ with $h,\rho\in\cO$ and $\deg_u(\rho)<\deg_u(\phi_f)$. On $\cO$ we have 
	$\deg_\tau=(j+d)\deg_u$, so $\phi_f(w-h)=\rho$ is possible only if $w=h$. Thus $w\in\cO$. 

	By Theorem \ref{thmA}, $D$ is the endomorphism algebra of a non-isotrivial Drinfeld module of rank $r'$ and 
	$A$-characteristic $(T)$ if and only if $r'=d(j+dn)$ with $n\geq1$. So $\phi$ has the smallest possible rank. If 
	$c\in\oK$ satisfies $c^{q^j-1}=x^{-1}$, then $c^{-1}\phi c$ reduces to $\psi_{dj}$ at $v$, and $dj=r-[D:F]$, in accordance 
	with Theorem \ref{thmAp}. For $d=3$ we get $D_3$ in rank $12$ from $u=x\tau+\tau^4$ and $D_3^{\mathrm{opp}}$ in rank $15$ from 
	$u=x\tau^2+\tau^5$. For $d=2$ and $j=1$, writing $T=y^2$, the module $\phi$ is the restriction to $A=\F_q[y^2]$ of the Drinfeld 
	$\F_q[y]$-module $y\mapsto x\tau+\tau^3$ of rank $3$. For $q=p$ this is the example of Pink \cite[Ex. 5.2 (c)]{PinkII}. 
\end{rem}

\begin{proof}[Proof of Theorem \ref{thmC}]
	Part (1) is Example \ref{ex1} and Remark \ref{remExOtherChar}, and part (2) is Example \ref{ex2}.
\end{proof}

\bibliographystyle{amsalpha}
\bibliography{EndAlgBib}

\end{document}